\documentclass[11pt]{article}
\usepackage{fullpage}

\usepackage{amssymb,amsmath,amsthm}
\usepackage{url}
\usepackage{verbatim}

\newtheorem{thm}{Theorem}
\newtheorem{cor}[thm]{Corollary}
\newtheorem{lem}[thm]{Lemma}

\theoremstyle{definition}

\usepackage{indentfirst}
\usepackage{xspace}
\usepackage[colorlinks=true,linkcolor=blue,citecolor=blue,urlcolor=blue]{hyperref}

\def\F{\mbox{\ensuremath{\mathcal F}}\xspace}
\def\G{\mbox{\ensuremath{\mathcal G}}\xspace}
\def\HH{\mbox{\ensuremath{\mathcal H}}\xspace}

\begin{document}
\title{\Large Note on polychromatic coloring of hereditary hypergraph families~III}
\author{D\"om\"ot\"or P\'alv\"olgyi
\thanks{ELTE E\"otv\"os Lor\'and University and Alfr\'ed R\'enyi Institute of Mathematics, Budapest, Hungary. Email: \texttt{domotor.palvolgyi@ttk.elte.hu}. Supported by the NRDI EXCELLENCE-24 grant no.~151504 Combinatorics and Geometry and by the ERC Advanced Grant no.~101054936 ERMiD.}}
\date{}
\maketitle

\noindent\textbf{Disclaimer.} This paper was written by ChatGPT with very limited feedback from me.

\begin{abstract}
For every integer $h\ge4$ we construct a hereditary family of hypergraphs with $m_2=h$ and
\[
m_3\ge 2(h-1)+2\left\lfloor\frac{h-1}{3}\right\rfloor+1
\approx \frac{8}{3}m_2.
\]
This improves the lower bound $m_3\ge 2m_2$ from the preceding notes for every $m_2\ge4$. We also show that the worst-case integrality gap of the weighted transversal problem over hereditary families with $m_2\le h$ is between $h-1$ and $5000h$. As a consequence, every such family admits strong $\varepsilon$-nets of size at most $5000h/\varepsilon$.
\end{abstract}


\noindent\textbf{Mathematics Subject Classification (MSC):} 05C15, 05C65

\medskip

On definitions and notation: In this note, all hypergraphs are finite. For $X\subseteq V(\HH)$, the trace $\HH[X]$ has edges $E\cap X$ for $E\in E(\HH)$. A subhypergraph of a trace is called a restricted subhypergraph. A hypergraph is $h$-heavy if each edge has size at least $h$, and it is $r$-uniform if each edge has size $r$. A coloring is polychromatic if every edge contains all colors. In particular, a polychromatic $2$-coloring is a proper $2$-coloring.

A family of hypergraphs $\F$ is hereditary if it contains every restricted subhypergraph of each of its members. Following \cite{Pal-note-I}, let $m_k(\F)$ be the least positive integer $m$ such that every $m$-heavy member of $\F$ has a polychromatic $k$-coloring, and set $m_k(\F)=\infty$ if no such $m$ exists. Thus $m_2(\F)$ is the threshold for proper $2$-colorability. For a single hypergraph $\HH$, we write $m_k(\HH)$ for the threshold of the hereditary family of all its restricted subhypergraphs.

The first note \cite{Pal-note-I} gives a $5$-uniform hypergraph with $m_2=3$ and $m_3=6$. The second note \cite{Pal-note-II} gives, for every $h\ge4$, a $(2h-1)$-uniform hypergraph with $m_2\le h$ and $m_3=2h$. We improve the latter bound for every $h\ge4$.

\begin{thm}\label{thm:construction}
For every integer $h\ge4$, put
\[
r(h)=2(h-1)+2\left\lfloor\frac{h-1}{3}\right\rfloor.
\]
There exists an $r(h)$-uniform hypergraph $\HH$ which has no polychromatic $3$-coloring but all of whose $h$-heavy restricted subhypergraphs are $2$-colorable. In particular, $m_2(\HH)\le h$ and $m_3(\HH)=r(h)+1$.
\end{thm}

To obtain the family in the abstract, add all $(h-1)$-edges on the vertex set of $\HH$ and take the hereditary family it generates. These added edges have no $h$-heavy traces, and the complete $(h-1)$-uniform subhypergraph on $2h-3$ vertices has no proper $2$-coloring. There are enough vertices since $r(h)\ge 2h$. The original hypergraph is still a subhypergraph, so $m_2=h$ and $m_3\ge r(h)+1$.

We now turn our attention to bounding the size of $\varepsilon$-nets.

For a hypergraph $\HH$ with nonempty edges and a nonnegative vertex cost $c$, let $\tau_c(\HH)$ be the minimum cost of a transversal, that is, a set meeting every edge. Its fractional relaxation is
\[
   \tau_c^*(\HH)=\min\left\{
      \sum_{v\in V(\HH)}c(v)x(v):
      x(v)\ge0,\ \sum_{v\in E}x(v)\ge1\ \text{for every }E\in E(\HH)
   \right\}.
\]
Let $G(h)$ be the supremum of $\tau_c(\HH)/\tau_c^*(\HH)$ over costs and hypergraphs with $\tau_c^*(\HH)>0$ and $m_2(\HH)\le h$.

\begin{thm}\label{thm:gap}
For every integer $h\ge2$,
\[
                        h-1\le G(h)\le 5000h .
\]
In particular, for every hereditary family $\F$ with $2\le m_2(\F)<\infty$, every $\HH\in\F$ with nonempty edges and every nonnegative vertex cost $c$ satisfy
\[
                \tau_c(\HH)\le 5000\,m_2(\F)\,\tau_c^*(\HH).
\]
\end{thm}

If $m_2(\F)\le h$, every member of $\F$ has VC-dimension at most $2(h-1)$. A shattered set of $2h-1$ vertices would give a restricted subhypergraph consisting of all $h$-subsets of the set, and every $2$-coloring would make one of those $h$-sets monochromatic. The Haussler--Welzl theorem \cite{HW87} thus gives strong $\varepsilon$-nets of size $O_h\bigl(\varepsilon^{-1}\log(1/\varepsilon)\bigr)$ as $\varepsilon\to0$. The following corollary removes the logarithmic factor.

\begin{cor}\label{cor:epsilon-nets}
Let $\F$ be a hereditary family with $m_2(\F)\le h<\infty$, where $h\ge2$. For every $\HH\in\F$, every probability measure $\mu$ on its vertices and every $0<\varepsilon\le1$, there is a strong $\varepsilon$-net of size at most $5000h/\varepsilon$. That is, there is a set of at most $5000h/\varepsilon$ vertices that meets every edge $E$ with $\mu(E)\ge\varepsilon$.
\end{cor}

\begin{proof}
Keep only the edges $E$ with $\mu(E)\ge\varepsilon$. The resulting subhypergraph belongs to $\F$ by heredity. The weights $x(v)=\mu(v)/\varepsilon$ form a fractional transversal of total weight $1/\varepsilon$. Apply Theorem~\ref{thm:gap} with unit costs.
\end{proof}

Linear growth of $m_k(\F)$ would also imply linear-size $\varepsilon$-nets for the uniform measure. Indeed, suppose $m_k(\F)\le Ck$ for all $k\ge2$. On an $n$-vertex member of $\F$, keep the edges of size at least $\varepsilon n$. If $\varepsilon n\ge2C$, put $k=\lfloor\varepsilon n/C\rfloor$. The retained subhypergraph is $m_k(\F)$-heavy, and a smallest color class of a polychromatic $k$-coloring meets every retained edge and has size at most $n/k\le2C/\varepsilon$. If $\varepsilon n<2C$, the entire vertex set satisfies the size bound. Thus the failure of linear-size $\varepsilon$-nets for the uniform measure, even of weak ones in a geometric realization, would have ruled out $m_k(\F)=O(k)$.

\section{The construction}\label{sec:construction}

We describe an auxiliary hypergraph $\G$ that will be used to define $\HH$. It was found by computer search to have the properties described in Lemma~\ref{lem:seed}. It has 14 vertices, five $2$-edges, and 140 $8$-edges. Its vertex set $V$ consists of two copies of $\mathbb Z_7$. The first has labels $i$ and the second has labels $i'$ for $i\in\mathbb Z_7$. Write $P_i=\{i,i'\}$. The short edges are $P_2,P_3,P_4,P_5,P_6$. To define the long edges, assign a pair $D(Q)\subset Q$ to each four-subset $Q\subset\mathbb Z_7$. The following five pairs $(Q,D(Q))$ specify the assignment by cyclic translation.
\[
\begin{array}{c|ccccc}
Q       &0123&0124&0134&0234&0135\\
D(Q)    &13&01&34&02&03
\end{array}
\]
Here a string such as $0123$ denotes the set $\{0,1,2,3\}$. First we assign all $35$ values of $D$ by translating the five table entries modulo $7$. For example, translating $Q=0123$ and $D(Q)=13$ by $6$ gives $Q=0126$ and the provisional value $D(0126)=02$. Once all $35$ values have been assigned, we change $D(0126)$ from $02$ to $06$ and $D(0124)$ from $01$ to $24$, and the other $33$ values are not changed.

For each $Q$, omit both points from each pair indexed by $D(Q)$ and one point from each pair indexed by $Q\setminus D(Q)$. This omits $2+2+1+1=6$ points. Each of the four choices of the two single points gives an $8$-edge by taking the complement in $V$. For example, when $Q=0126$, we omit all of $P_0$ and $P_6$. If we also omit $1$ and $2$, the resulting edge is
\[
P_3\cup P_4\cup P_5\cup\{1',2'\}.
\]
Each such edge determines $Q$, which consists exactly of the indices of the four pairs that are not fully contained in the edge. For a fixed $Q$, the four choices give distinct edges. Thus there are $35\cdot4=140$ long edges.

\begin{lem}\label{lem:seed}
The hypergraph $\G$ has the following properties.
\begin{enumerate}
\item Every subset of at most four points is disjoint from some long edge and from some short edge.
\item For every $X\subseteq V$, there is a red-blue coloring of $X$ in which every monochromatic trace of a long edge has size at most three and every monochromatic trace of a short edge has size at most one.
\end{enumerate}
\end{lem}

\begin{proof}
The five short edges are pairwise disjoint, so no four points meet all of them. For the long edges it suffices to test four-sets, since a smaller set can be extended to a four-set. The deterministic program in Appendix~\ref{app:code} checks all $\binom{14}{4}=1001$ four-sets against the $140$ long edges. For the second assertion, it checks all $2^{14}=16384$ choices of $X$. For each $X$ it looks for a coloring that splits every long trace of size at least four and every short edge contained in $X$.
\end{proof}

\begin{proof}[Proof of Theorem~\ref{thm:construction}]
Write $h-1=3q+t$, where $q=\lfloor(h-1)/3\rfloor$ and $t=(h-1)\bmod 3$, so $0\le t\le2$. Take $3(q+t-1)+1$ disjoint copies of $\G$. An edge of $\HH$ is the union of one long edge from each of $q$ copies and one short edge from each of $t$ other copies. Every such edge has size $8q+2t=2(h-1)+2q=r(h)$.

Assume that there is a polychromatic $3$-coloring. In each copy of $\G$, choose a smallest color class, which has at most four vertices. By the pigeonhole principle, the same color is chosen in at least $q+t$ copies. By Lemma~\ref{lem:seed}, we can choose a long edge avoiding that color in each of $q$ of those copies and a short edge avoiding it in each of $t$ others. Their union is an edge of $\HH$ that misses the chosen color, a contradiction.

Now let $X$ be any subset of the vertex set. Color its intersection with each copy as in Lemma~\ref{lem:seed}. A monochromatic trace of an edge of $\HH$ has at most three points from each long edge and at most one from each short edge. Its size is at most $3q+t=h-1$. Thus $m_2(\HH)\le h$. The absence of a polychromatic $3$-coloring gives $m_3(\HH)>r(h)$, and no restricted subhypergraph has an edge larger than $r(h)$, so $m_3(\HH)=r(h)+1$.
\end{proof}

\section{The weighted transversal gap}\label{sec:gap}

For the lower bound in Theorem~\ref{thm:gap}, take all $(h-1)$-sets on $n$ vertices, with unit costs. The condition on $h$-heavy restrictions is vacuous. The transversal number is $n-h+2$ and the fractional transversal number is $n/(h-1)$, by symmetry or by averaging the edge constraints. Their ratio tends to $h-1$ as $n$ tends to infinity.

For the upper bound we first treat costs that are integer powers of two. In the sampling argument, each cost class is sampled separately. We bound the cost of the sampling errors simultaneously for all red-blue labelings, so that we can apply the bound to a coloring supplied by the hereditary hypothesis.

The argument is related to the double-sampling proof of the $\varepsilon$-net theorem of Haussler and Welzl \cite{HW87}, also presented by Matou\v{s}ek \cite[Chapter~10]{Mat02}. A related use of concentration and a union bound over patterns induced by LP dual solutions appears in Agrawal, Wang, and Ye \cite[Lemma~2]{AWY14}. The step specific to our proof is to minimize the objective defining $F_{i+1}$ over the larger feasible set defining $F_i$. This fixes a dual packing before revealing the draws in cost class $i$ whose positions are labeled blue and gives a comparison that telescopes over the cost classes.

\begin{proof}[Proof of the upper bound in Theorem~\ref{thm:gap}]
We may delete the zero-cost vertices and all edges meeting them, since one can include these vertices in a transversal at no cost. This preserves both optima and the hereditary hypothesis. Round each remaining cost up to a power of two. The integral optimum can only increase, and the fractional optimum grows by at most a factor of two. It thus suffices to prove the bound $2500h$ for costs that are powers of two.

For the given costs, pass to a subfamily of edges maximizing the ratio of the two optima among nonempty subfamilies. Denote the ratio by $R$. Fix a rational optimal fractional transversal $x$ and write $C=\tau_c^*(\HH)=\sum_{v\in V(\HH)}c(v)x(v)$. We may take sufficiently many disjoint copies of this instance and of $x$. Keep the notation $x$ and $C$ for the fractional transversal and its cost on the union. The ratio remains $R$ and the same upper bound $R$ applies to every nonempty subfamily, by additivity of both optima over the copies. Thus we can assume that $10h\sum_{v:c(v)=a}x(v)$ is an integer at least $1000$ for every cost $a$ with positive $x$-mass. The hereditary hypothesis also persists under disjoint union.

List these cost classes in increasing order as $c_1,\ldots,c_k$. For each $i$, draw
\[
                  n_i=10h\sum_{v:c(v)=c_i}x(v)
\]
positions from that class with replacement, with vertex probabilities proportional to $x(v)$. All draws are independent. The total cost of the positions, counted with repetitions, is $\sum_i c_i n_i=10hC$. For any outcome, its positions admit a red-blue labeling in which each edge with at least $h$ sampled occurrences has both colors. To see this, first split the occurrences of each vertex sampled at least twice between the colors. Apply the hereditary hypothesis on the set of vertices sampled exactly once to color those positions.

A feasible dual packing $y$ assigns a nonnegative weight $y_E$ to every edge and has load
\[
           \ell_y(v)=\sum_{E\ni v} y_E\le c(v).
\]
Its total weight is at most $C$ by LP duality. For a labeling of the sampled positions define the slacks
\[
 d_i(y)=\sum_{v:c(v)=c_i}x(v)\bigl(c_i-\ell_y(v)\bigr),\qquad
 e_i(y)=\sum_{\text{blue positions }v\text{ in class }i}
                   \bigl(c_i-\ell_y(v)\bigr).
\]
Let $Y_i$ be the set of feasible dual packings supported on edges avoiding all sampled positions in classes $j<i$ and all red positions in classes $j\ge i$. In particular $Y_{k+1}$ avoids every sampled position. Put
\[
 N_i=\sum_{j\ge i}n_j,\qquad
 F_i=\min_{y\in Y_i}\left(30\cdot10h\sum_{j<i}d_j(y)+\sum_{j\ge i}e_j(y)\right).
\]

We claim that there is a sample for which, simultaneously for every labeling and every $i\le k$,
\begin{equation}\label{eq:step}
                     F_i\le F_{i+1}+0.14c_i N_i .
\end{equation}
Fix $i$ and fix the labels in classes $j\ge i$. Condition on the draws outside class $i$ and the draws in class $i$ whose positions are labeled red. Choose $\widehat y$ to minimize
\[
                     30\cdot10h\sum_{j\le i}d_j(y)+\sum_{j>i}e_j(y)
\]
over $Y_i$. This objective and $Y_i$ are now fixed, independently of the draws in class $i$ whose positions are labeled blue. Since $Y_{i+1}\subseteq Y_i$, the value at $\widehat y$ is at most $F_{i+1}$.

For a draw at $v$ in a position labeled blue in class $i$, its normalized slack $Z=(c_i-\ell_{\widehat y}(v))/c_i$ lies in $[0,1]$. By convexity, $\exp(5Z)\le 1+(e^5-1)Z$. The conditional expectation of the sum of the normalized slacks over the positions labeled blue is at most $10h\,d_i(\widehat y)/c_i$, since the number of positions labeled blue is at most $n_i$. Independence and $1+z\le e^z$ give
\[
  \mathbb E\!\left[\exp\!\left(5e_i(\widehat y)/c_i\right)\right]
          \le \exp\!\left((e^5-1)\cdot10h\,d_i(\widehat y)/c_i\right).
\]
Here and in the following probability bound the conditioning is as above. Since $e^5-1<5\cdot30$, Markov's inequality gives
\[
 \Pr\!\left(e_i(\widehat y)>30\cdot10h\,d_i(\widehat y)+0.14c_iN_i\right)
       \le e^{-0.7N_i}.
\]
Outside this event we can evaluate the objective of $F_i$ at $\widehat y$ and obtain \eqref{eq:step}. The bound also holds without conditioning.

Only the labels in classes $j\ge i$ enter the two minima, giving $2^{N_i}$ choices for the union bound at step $i$. Since $n_i\ge1000$, we have $N_i\ge1000(k-i+1)$. Thus the total failure probability is at most
\[
   \sum_{i=1}^k (2e^{-0.7})^{N_i}
       \le \sum_{j\ge1}(2e^{-0.7})^{1000j}<1 .
\]
Fix an outcome for which \eqref{eq:step} holds for all labelings.

Use the hereditary red-blue labeling of the sampled positions described above. Interchange the colors if necessary so that the blue positions have cost at least $10hC/2$. An edge supporting a packing in $Y_1$ avoids red and therefore has at most $h-1$ blue positions. Consequently
\[
 F_1\ge \frac{10hC}{2}-(h-1)C>0.4\cdot10hC.
\]
For the other endpoint, discard all edges hit by the sample. Adding the sampled vertices to a transversal of the remaining subfamily gives a transversal of $\HH$ at additional cost at most $10hC$. The remaining subfamily therefore has integral optimum at least $RC-10hC$ and fractional optimum at least $C-10hC/R$, by the choice of $R$. The last inequality also holds for an empty subfamily, in which case $RC\le 10hC$. By duality there is a packing $y\in Y_{k+1}$ of weight at least $C-10hC/R$. Since $x$ is feasible,
\[
 \sum_i d_i(y)=C-\sum_E y_E\sum_{v\in E}x(v)
             \le C-\sum_E y_E\le 10hC/R .
\]
It follows that $F_{k+1}\le 30(10h)^2C/R$.

Finally, the different costs are powers of two, so $\sum_{i\le j}c_i\le2c_j$. Summing \eqref{eq:step} and using
\[
  \sum_i c_iN_i=\sum_j n_j\sum_{i\le j}c_i\le 2\cdot10hC
\]
yields
\[
              0.4\cdot10hC<F_1\le 30(10h)^2C/R+0.28\cdot10hC.
\]
Therefore $R<30(10h)/0.12=2500h$, as required.
\end{proof}

Now we state one more strengthening of Theorem~\ref{thm:gap}. The proof is a straightforward adaptation of the argument above, so we omit the details. For $k\ge2$, let $p_k(\F)$ be the least positive integer $h$ such that every $h$-heavy member of $\F$ admits a proper $k$-coloring, in which no edge is monochromatic. Set $p_k(\F)=\infty$ if no such $h$ exists. Thus $p_2(\F)=m_2(\F)$.

\begin{thm}\label{thm:proper}
For each $k\ge2$ there is a constant $C_k$ such that, if $\F$ is hereditary and $p_k(\F)\le h<\infty$, then
\[
  \tau_c(\HH)\le C_k h\,\tau_c^*(\HH)
\]
for every member $\HH$ with nonempty edges and every nonnegative vertex cost $c$. For every probability measure on the vertices of a member of $\F$ and every $0<\varepsilon\le1$, there is a strong $\varepsilon$-net of size at most $C_k h/\varepsilon$.
\end{thm}

The same shattering argument gives a VC-dimension bound of $k(h-1)$ whenever $p_k(\F)\le h$.

\appendix
\section{Verification code}\label{app:code}

The following Python program verifies Lemma~\ref{lem:seed}. It represents vertex sets by bit masks and uses no external packages.

\begin{scriptsize}
\setlength{\baselineskip}{7pt}
\begin{verbatim}
from itertools import combinations, product
def mask(indices):
    return sum(1 << int(i) for i in indices)
def shift(indices, j):
    return sum(1 << ((int(i) + j) % 7) for i in indices)
D = {}
for Q, pair in [("0123", "13"), ("0124", "01"), ("0134", "34"),
                ("0234", "02"), ("0135", "03")]:
    for j in range(7):
        key = shift(Q, j)
        assert key not in D
        D[key] = shift(pair, j)
assert len(D) == 35
D[mask("0126")] = mask("06")
D[mask("0124")] = mask("24")
V = (1 << 14) - 1
long_edges = set()
for Q, pair in D.items():
    assert Q.bit_count() == 4 and pair.bit_count() == 2
    assert Q & pair == pair
    i, j = [k for k in range(7) if (Q ^ pair) & (1 << k)]
    omit = sum(3 << (2*k) for k in range(7) if pair & (1 << k))
    for a, b in product(range(2), repeat=2):
        long_edges.add(V ^ (omit | (1 << (2*i+a)) | (1 << (2*j+b))))
short_edges = [3 << (2*k) for k in range(2, 7)]
assert len(long_edges) == 140
assert all(e.bit_count() == 8 for e in long_edges)
for four in combinations(range(14), 4):
    assert any(e & mask(four) == 0 for e in long_edges)
def admissible(S):
    heavy = {e & S for e in long_edges if (e & S).bit_count() >= 4}
    heavy.update(e for e in short_edges if e & S == e)
    if not heavy:
        return True
    R = S
    while R:
        B = S ^ R
        if all(e & R and e & B for e in heavy):
            return True
        R = (R - 1) & S
    return False
for S in range(1 << 14):
    assert admissible(S), S
print("140 long edges, 1001 four-sets, 16384 restrictions verified.")
\end{verbatim}
\end{scriptsize}

\end{document}